\documentclass[a4paper,fleqn]{cas-dc}

\usepackage[authoryear,longnamesfirst]{natbib}

\newdefinition{remark}{Remark}
\usepackage[nameinlink,capitalise,noabbrev]{cleveref} 
\usepackage{enumitem}
\usepackage{amsthm}

\newcommand{\RR}{\ensuremath{\mathbb{R}}}

\newcommand{\AAA}{\ensuremath{\mathcal{A}}}

\newcommand{\GGG}{\ensuremath{\mathcal{G}}}

\newcommand{\XXX}{\ensuremath{\mathcal{X}}}

\newcommand{\IIII}{\ensuremath{\mathfrak{I}}}
\newcommand{\LLLL}{\ensuremath{\mathfrak{L}}}

\crefformat{equation}{(#2#1#3)}
\crefformat{enumi}{#2#1#3}
\crefmultiformat{enumi}{#2#1#3}{ and ~#2#1#3}{, #2#1#3}{ and~#2#1#3}
\crefname{problem}{problem}{problems}
\Crefname{problem}{Problem}{Problems}

\crefformat{figure}{#2figure~#1#3}
\Crefformat{figure}{#2Figure~#1#3}

\newcommand{\lccref}[1]{\hyperref[{#1}]{\lcnamecref{#1}~\labelcref{#1}}} 

\newcommand\lb[2]{\left[#1,#2\right]}              
\renewcommand\vec[1]{\frac{\partial}{\partial #1}} 
\newcommand\distrib[1]{\mathrm{span}\left\{#1\right\}} 
\newcommand\ad[2]{\mathrm{ad}_{#1}#2}                
\newcommand\adk[3]{\mathrm{ad}_{#1}^{#2}#3}          
\newcommand\vectR[1]{\mathrm{vect}_{\RR}\left\{#1\right\}}       
\newcommand\id{\mathrm{Id}}       					 

\newcommand\diff{\mathrm{d}} 

\newcommand\Sm{\ensuremath{S}}

\makeatletter
\newcommand*\owedge{\mathpalette\@owedge\relax}
\newcommand*\@owedge[1]{%
  \mathbin{%
    \ooalign{%
      $#1\m@th\bigcirc$\cr
      \hidewidth$#1\m@th\wedge$\hidewidth\cr
    }%
  }%
}
\makeatother

\newcommand{\tred}[1]{\textcolor{red}{#1}}
\renewcommand{\tred}{}

\newtheorem{theorem}{Theorem}

\newtheorem{proposition}{Proposition}
\newtheorem{lemma}{Lemma}

\let\OldItem\item
\newcommand{\MyItem}[2][]{}%
\newenvironment{myenumerate}[1][]{%
    \renewcommand{\item}[2][]{%
        \begin{enumerate}[#1,label={##1},ref={##1}]%
            \OldItem {##2}%
        \end{enumerate}%
    }%
}{%
}%

\newcommand{\RomanNumeralCaps}[1]
    {\MakeUppercase{\romannumeral #1}}

\graphicspath{{img/}}

\begin{document}


\shorttitle{Symmetries of conic systems}    

\shortauthors{T. Schmoderer \& W. Respondek}  

\title[mode = title]{Null-forms of conic systems in $\RR^3$ are determined by their symmetries}  



%

\author[1]{Timothée Schmoderer}[bioid=1, orcid=0000-0003-1748-1763]

\cormark[1]


\ead{timothee.schmoderer@insa-rouen.fr}

\ead[url]{tschmoderer.github.io}


\affiliation[1]{organization={Laboratoire de Mathématiques},
            addressline={INSA Rouen Normandie}, 
            city={Rouen},
            country={France}}

\author[1]{Witold Respondek}[bioid=2]


\ead{witold.respondek@insa-rouen.fr}




\cortext[1]{Corresponding author}



\begin{abstract}
We address the problem of characterisation of null-forms of conic $3$-dimensional systems, that is, control-affine systems whose field of admissible velocities forms a conic (without parameters) in the tangent space. Those systems have been previously identified as the simplest control systems under a conic nonholonomic constraint or as systems of zero curvature. In this work, we propose a direct characterisation of null-forms of conic systems among all control-affine systems by studying the Lie algebra of infinitesimal symmetries. Namely, we show that the Lie algebra of infinitesimal symmetries characterises uniquely null-forms of conic systems. 
\end{abstract}




\maketitle
\section{Introduction}\label{sec:intro}
In \cite{schmoderer2021Conicnonholonomicconstraints}, we propose a characterisation and a classification of single input control-affine systems which admit conic sets as admissible velocities. In particular, we give a characterisation of all $3$-dimensional control-affine systems with scalar control that are feedback equivalent, locally around $(x_0,y_0,w_0)\in\RR^3$, to one of the following conic null-forms 
\begin{align*}
    \Sigma_E:&\left\{\begin{array}{cl}
        \dot{x} &= \cos(w) \\
        \dot{y} &= \sin(w) \\ 
        \dot{w} &= u
    \end{array}\right., \quad
    \Sigma_H:\left\{\begin{array}{cl}
        \dot{x} &= \cosh(w) \\
        \dot{y} &= \sinh(w) \\ 
        \dot{w} &= u
    \end{array}\right.,\\
    &\textrm{and }\quad
    \Sigma_P:\left\{\begin{array}{cl}
        \dot{x} &= w^2 \\
        \dot{y} &= w \\ 
        \dot{w} &= u
    \end{array}\right..
\end{align*}
\noindent
    {The parabolic system $\Sigma_P$ actually depends on the nature of the point $(x_0,y_0,w_0)$ around which the systems $\Sigma_P$ is considered: an equilibrium or not. As a consequence, $\Sigma_P$ has two distinct normal forms, for $w_0=0$ and $w_0\neq0$, and we will specify this distinction just before our main result \cref{thm:m1_symmetries_characterize_systems}.}

We call the systems $\Sigma_E$, $\Sigma_H$, and $\Sigma_P$, \emph{null-forms} because of the absence of any parameters, functional or real-continuous, in their expression (in particular, observe that those system are trivial in the sense of \cite{serres2009Controlsystemszero}). We call the systems $\Sigma_E$, $\Sigma_H$, and $\Sigma_P$, elliptic, hyperbolic, and parabolic, respectively, because {in the manifold $\RR^2$, equipped with local coordinates $(x,y)$, the trajectories $(x(t),y(t))$ of $\Sigma_E$, $\Sigma_H$, and $\Sigma_P$ satisfy a nonholonomic constraint of the form 
\begin{align*}
    \Sm_E=\dot{x}^2+\dot{y}^2-1&=0,\quad \Sm_H=\dot{x}^2-\dot{y}^2-1=0,\\
    \textrm{and}&\quad\Sm_P=\dot{y}^2-\dot{x}=0,
\end{align*}
\noindent
respectively. \tred{The importance of the null-forms $\Sigma_E$, $\Sigma_H$, and $\Sigma_P$, is multiple. First}, $\Sm_E$, resp. $\Sm_H$, and resp. $\Sm_P$, describes ellipses, resp. hyperbolas, and resp. parabolas, in the tangent bundle $T\RR^2$\tred{, that is, the simplest non-holonomic constraints that are not affine with respect to velocities}; 
see \cite{schmoderer2021Conicnonholonomicconstraints} for a detailed analysis of the link between quadratic equations on $T\RR^2$ and control systems.} Second, the elliptic and hyperbolic systems $\Sigma_E$ and $\Sigma_H$ show up in the seminal trend-setting paper \cite{agrachev1998FeedbackInvariantOptimalControl} as \tred{simple }control-affine systems of zero curvature \tred{(a crucial feedback invariant of control-affine system on $3$-dimensional manifolds with scalar control)}. Third, $\Sigma_E$ describes the well-known Dubins model of a car \cite{dubins1957CurvesMinimalLength}, whereas $\Sigma_H$ is its hyperbolic counterpart \cite{monroy-perez1998NonEuclideanDubinsProblem}, and while the parabolic system $\Sigma_P$ is studied in \cite{schmoderer2021Conicnonholonomicconstraints}. \tred{Fourth, last but not least, the null-forms $\Sigma_E$, $\Sigma_H$, and $\Sigma_P$, have the simplest possible Lie ideals $\mathfrak{L}_0$ among all single-input control-affine systems with $3$-dimensional state space. To justify that statement, recall that for the system $\dot{x}=f(x)+g(x)u$ we define its Lie algebra $\mathfrak{L}=\left\{f,g\right\}_{\mathrm{LA}}$ and its Lie ideal $\LLLL_0$ (generated by $g$ in $\LLLL$), the former being responsible for accessibility and the latter for strong accessibility, see e.g. \cite{isidori1995NonlinearControlSystems} and \cite{nijmeijer1990NonlinearDynamicalControl}. For $\Sigma_E$, $\Sigma_H$, and $\Sigma_P$, the Lie ideal is the simplest possible, namely, it is the $3$-dimensional Lie algebra $\LLLL_0=\vectR{g,\ad{f}{g},\lb{g}{\ad{f}{g}}}$ such that $\IIII=\vectR{\ad{f}{g},\lb{g}{\ad{f}{g}}}$ is an abelian ideal of $\LLLL_0$ and the linear map $\ad{g}{}:\IIII\rightarrow\IIII$ has eigenvalues $\pm i$ for $\Sigma_E$, $\pm1$ for $\Sigma_H$, or a double zero eigenvalue (and is nilpotent) for $\Sigma_P$; compare \cref{rm:lie_ideal}.  } 

We will denote by $\Sigma_Q$ the set $\{\Sigma_E,\Sigma_H,\Sigma_P\}$ of the three above null-forms conic systems when we discuss properties that do not require to distinguish a specific system. \tred{Consider the system $\Xi_Q:\dot{x}=f_Q(x,w)$, where $x\in\RR^2$, $w\in\RR$, and $f_Q$ takes one of the three forms $f_E=A(x)\cos(w)+B(x)\sin(w)+C(x)$, or $f_H=A(x)\cosh(w)+B(x)\sinh(w)+C(x)$, or $f_P=A(x)w^2+B(x)w+C(x)$, where $A$, $B$, and $C$ are smooth vector fields on $\RR^2$. Observe, for the null-forms $\Sigma_Q$, that interpreting $w\in\RR$ as a scalar control that enters non-linearly, the three forms $\Sigma_Q$ can be considered as $\Xi_Q$ for which the vector fields $A$ and $B$ are simultaneously rectified and $C$ is annihilated. In \cite{schmoderer2021Conicnonholonomicconstraints}, see also \cite{schmoderer2018Studycontrolsystems}, we described (by identifying suitable feedback invariants) the systems $\Xi_Q$ with $A$, $B$, and $C$ of the just mentioned form, thus we characterised the null-forms $\Sigma_Q$. The resulst of the present paper differ from those of \cite{schmoderer2021Conicnonholonomicconstraints} in two aspects. First, in \cite{schmoderer2021Conicnonholonomicconstraints}, we characterised $\Sigma_Q$ among all conic systems of the form $\Xi_Q$, while in the present paper, we propose a characterisation that directly applies to any   } 
control-affine system (on a $3$D-manifold with scalar control). Second, and more important, the methodology that we apply to describe 
our null-forms is based on the study of the Lie algebra of infinitesimal symmetries (see below for a precise definition), which turns out to uniquely determine all null-forms $\Sigma_Q$ (with additional regularity conditions for $\Sigma_P$). \tred{To summarise, the three null forms $\Sigma_Q$ appeared in our previous work \cite{schmoderer2021Conicnonholonomicconstraints}, where we proposed their characterisation as subclasses of control-affine systems satisfying a quadratic nonholonomic constraint. In particular, the null-forms $\Sigma_Q$ show up as the simplest models, i.e. with no functional nor real parameters, in our classification and therefore we expected that they are the most symmetric. In the present paper, we propose their study via the infinitesimal symmetries and prove that their symmetries (which, indeed, form basic $3$-dimensional Lie algebras) completely determine the forms $\Sigma_Q$ up to feedback transformations.}
{
\subsection{Outline of the paper}
The paper is organised as follows. In the next subsection, we will introduce all notions of control theory that we will need for our results. In particular, we give a definition of feedback equivalence and of infinitesimal symmetries. Next,  we will first compute the Lie algebra of infinitesimal symmetries of null-forms of conic systems $\Sigma_Q$ and identify those Lie algebras in the well-known Bianchi classification of $3$-dimensional Lie algebras. Afterwards, we will enunciate and prove our main result, which shows that the Lie algebras of infinitesimal symmetries of null-forms of conic systems completely determine corresponding systems. 
}
\subsection{Preliminaries}
Throughout the paper, the word "smooth" will always mean $C^{\infty}$-smooth, and all objects (vector fields, differential forms, functions, manifolds) are assumed to be smooth. 
\paragraph{Control-affine systems.}
We consider control-affine systems $\Sigma$ of the form
\begin{align*}
    \Sigma\,:\,\dot{\xi}=f(\xi)+g(\xi)u,\quad u\in\RR,
\end{align*}
\noindent
where the state $\xi$ belongs to a smooth $3$-dimensional manifold $\XXX$ (or an open subset of $\RR^3$ since all our results are local), and $f$ and $g$ are smooth vector fields on $\XXX$ (smooth sections of the tangent bundle $T\XXX$). We denote a control-affine system by the pair $\Sigma=(f,g)$, we set $\GGG=\distrib{g}$, the distribution spanned by the vector field $g$, and we set $\GGG^1=\distrib{g,\lb{f}{g}}$. 
We call two control affine systems $\Sigma=(f,g)$ and $\tilde{\Sigma}=(\tilde{f},\tilde{g})$ on $\tilde{\XXX}$ feedback equivalent, if there exists a diffeomorphism $\phi :\XXX\rightarrow\tilde{\XXX}$ and smooth functions $\alpha(\xi)$ and $\beta(\xi)$, satisfying $\beta(\xi)\neq0$, and such that 
\begin{align*}
    \tilde{f}=\phi_*(f+g\alpha)\quad\textrm{and}\quad\tilde{g}=\phi_*(g\beta),
\end{align*}
\noindent
where $\phi_*$ denotes the tangent map of $\phi$, that is $(\phi_*f)(\tilde{\xi})=\frac{\partial \phi}{\partial \xi}\left(\phi^{-1}(\tilde{\xi})\right)\cdot f\left(\phi^{-1}(\tilde{\xi})\right)$. If $\phi$ is defined locally around $\xi_0$ and $\tilde{\xi}_0=\phi(\xi_0)$, then we say that $\Sigma$ and $\tilde{\Sigma}$ are locally feedback equivalent at $\xi_0$ and $\tilde{\xi}_0$, respectively. A feedback transformation will be denoted by the triple $(\phi,\alpha,\beta)$ or simply $(\alpha,\beta)$ if $\phi=\id$. 

%
%
\paragraph{Infinitesimal symmetries.}
We introduce the notion of symmetries in the case of single input control-affine system (see \cite{respondek2002Nonlinearizablesingleinputcontrol,grizzle1985structurenonlinearcontrol} for a detailed introduction). For a control-affine system $\Sigma=(f,g)$, we define the field of admissible velocities $\AAA$ as
\begin{align*}
    \AAA(\xi)=\{f(\xi)+g(\xi)u\, :\, u\in\RR\}\subset T_{\xi}\XXX .
\end{align*}
\noindent
We say that a diffeomorphism $\phi\, :\, \XXX\rightarrow\XXX$ is a \emph{symmetry} of $\Sigma$ if it preserves the field of affine lines $\AAA$ (equivalently, the affine distribution $\AAA=f+\GGG$), that is,
\begin{align*}
    \phi_*\AAA=\AAA.
\end{align*}
\noindent
We say that a vector field $v$ on $\XXX$ is an infinitesimal symmetry of $\Sigma=(f,g)$ if the (local) flow $\gamma^v_t$ of $v$ is a (local) symmetry, for any $t$ for which it exists, that is \tred{the (local) diffemorphism $\gamma_t^v$ satisfies} $(\gamma_t^v)_*\AAA=\AAA$. Consider the system $\Sigma=(f,g)$ and recall that $\GGG$ is the distribution spanned by the vector field $g$. Since $(\gamma_t^v)_*(f+\GGG)=f+\GGG$, it follows that $(\gamma_t^v)_*$ preserves $\GGG$ and preserves $f$ modulo $\GGG$, which immediately gives the following characterisation of infinitesimal symmetries. 
\begin{proposition}\label{prop:charac_lie_algebra_infitesimal_symmetries}
A vector field $v$ is an infinitesimal symmetry of the control-affine system $\Sigma=(f,g)$ if and only if 
\begin{align*}
    \lb{v}{g}=0\hspace{-0.5em}\mod\GGG\quad\textrm{and}\quad\lb{v}{f}=0\hspace{-0.5em}\mod\GGG.
\end{align*}
\end{proposition}
\noindent
By the Jacobi identity, it is easy to see that if $v_1$ and $v_2$ are infinitesimal symmetries, then so is $\lb{v_1}{v_2}$, hence the set of all infinitesimal symmetries forms a real Lie algebra. Notice that the Lie algebra of infinitesimal symmetries is attached to the affine distribution $\AAA=f+\GGG$ and not to {a particular pair $(f,g)$. Different pairs $(f,g)$ related via feedback transformations $(\alpha,\beta)$ define the same $\AAA$ and thus give rise to the same Lie algebra of infinitesimal symmetries which, therefore, is a feedback invariant object attached to~$\Sigma$.} 
%
%

\section{Main results}
In this section, we prove that the Lie algebras of infinitesimal symmetries of the systems $\Sigma_E$, $\Sigma_H$, and $\Sigma_P$, determine the corresponding systems among all control-affine systems. With the help of \cref{prop:charac_lie_algebra_infitesimal_symmetries}, we give in the next lemma the Lie algebras of infinitesimal symmetries of the systems $\Sigma_E$, $\Sigma_H$, and $\Sigma_P$.
\begin{lemma}[Lie algebra of infinitesimal symmetries of $\Sigma_Q$]\label{lem:lie_algebra_symmetries}
The Lie algebra of infinitesimal symmetries of $\Sigma_E$, $\Sigma_H$, and $\Sigma_P$ are respectively:
\begin{align*}
    \mathfrak{L}_E &= \vectR{\vec{x},\vec{y},y\vec{x}-x\vec{y}-\vec{w}},\\
    \mathfrak{L}_H &= \vectR{\vec{x},\vec{y},y\vec{x}+x\vec{y}+\vec{w}},\\
        \mathfrak{L}_P &= \vectR{\vec{x},\vec{y},2x\vec{x}+y\vec{y}+w\vec{w}}.
\end{align*}
\end{lemma}

We will denote by $\mathfrak{L}_Q$ the set of the three above Lie algebras in order to state general facts about them. Those three Lie algebras share the property of having the abelian Lie ideal $\vectR{\vec{x},\vec{y}}$, which corresponds to the fact that the systems $\Sigma_Q$ are invariant under translations $(x,y)\mapsto(x+c,y+d)$, with $c,d\in\RR$. The third vector field plays the role of an Euler vector field for $\mathfrak{L}_P$ and is an infinitesimal rotation (\tred{elliptic} or hyperbolic) in the case of $\mathfrak{L}_E$ and $\mathfrak{L}_H$.

\begin{proof}
    The proof is a direct calculation using the characterisation of infinitesimal symmetries given in \cref{prop:charac_lie_algebra_infitesimal_symmetries}.
    Let $v=v_1\vec{x}+v_2\vec{y}+v_3\vec{w}$, where $v_i=v_i(x,y,w)$, be an infinitesimal symmetries of $\Sigma_Q$. Then, using that $v$ is a symmetry of $\GGG$, we deduce that $v_1=v_1(x,y)$ and $v_2=v_2(x,y)$. We now distinguish the three conic systems. 
    
    For $\Sigma_E$, we have $f=\cos(w)\vec{x}+\sin(w)\vec{y}$, hence $\lb{v}{f}\in\GGG$ implies 
    \begin{align*}
        \cos(w)\frac{\partial v_1}{\partial x}+\sin(w)\frac{\partial v_1}{\partial y}&=-\sin(w)v_3,\\ 
        \cos(w)\frac{\partial v_2}{\partial x}+\sin(w)\frac{\partial v_2}{\partial y}&=\cos(w)v_3.
    \end{align*}
    \noindent
    Multiplying the first equation by $\cos(w)$ and inserting the second equation, we obtain
    \begin{multline*}
        \left(\frac{\partial v_1}{\partial x}-\frac{\partial v_2}{\partial y}\right)\cos(2w)+\left(\frac{\partial v_1}{\partial y}+\frac{\partial v_2}{\partial x}\right)\sin(2w)\\+\left(\frac{\partial v_1}{\partial x}+\frac{\partial v_2}{\partial y}\right)=0,
    \end{multline*}
    \noindent
    which yields, first, $v_1=v_1(y)$ and $v_2=v_2(x)$, and, second, $v_1(y)=ay+b$ and $v_2(x)=-ax+c$, where $a,b,c\in\RR$. Finally, taking $(a,b,c)=(1,0,0)$, $(a,b,c)=(0,1,0)$, and $(a,b,c)=(0,0,1)$ we obtain $\LLLL_E$.
    
    For $\Sigma_H$, it is an analogous calculation as above but with hyperbolic functions.
    
    For $\Sigma_P$, we have $f=w^2\vec{x} +w\vec{y}$, hence $\lb{v}{f}\in\GGG$ implies 
    \begin{align*}
        w^2\frac{\partial v_1}{\partial x}+w\frac{\partial v_1}{\partial y}&=2wv_3,\\
        w^2\frac{\partial v_2}{\partial x}+w\frac{\partial v_2}{\partial y}&=v_3.
    \end{align*}
    \noindent
    Inserting the second equation in the first, we obtain a polynomial of degree $3$ in $w$ that identically vanishes. Thus each of its coefficient vanishes, implying $\frac{\partial v_2}{\partial x}=0$, $\frac{\partial v_1}{\partial y}=0$, and $\frac{\partial v_1}{\partial x}=2\frac{\partial v_2}{\partial y}$. Therefore, we conclude that $v=(2ax+b)\vec{x}+(ay+c)\vec{y}+aw\vec{w}$, where $a,b,c\in\RR$. Finally, taking $(a,b,c)=(1,0,0)$, $(a,b,c)=(0,1,0)$, and $(a,b,c)=(0,0,1)$ we obtain $\LLLL_P$. 
\end{proof}

Based on the commutativity relations of the generators of $\mathfrak{L}_Q$, we identify the three Lie algebras in the classification of $3$-dimensional Lie algebras as presented by Winternitz in \cite{bowers2005Classificationthreedimensionalreal}. First, $\mathfrak{L}_E\cong L(3,4,0)$ is the Lie algebra of the euclidean group $\textrm{E}(2)$, that is, the group of affine transformations of $\RR^2$ which preserve the euclidean metric. {In other words, $\LLLL_E$ is $\mathrm{\RomanNumeralCaps{7}}_0$ in the \cite{bianchi1898Suglispaziitre} classification of $3$-dimensional Lie algebras.} Second, $\mathfrak{L}_H\cong L(3,2,-1)$ is the Lie algebra of the Poincar\'e group $\textrm{P}(1,1)$, of affine transformations of $\RR^2$ which preserve the Lorentz metric{, i.e. $\LLLL_H$ is $\mathrm{\RomanNumeralCaps{6}}_0$ is the Bianchi classification.} Finally, $\mathfrak{L}_P\cong L(3,2,2)$ with no immediate interpretation {and $\LLLL_P$ is one of $\mathrm{\RomanNumeralCaps{6}}$-algebras of the Bianchi classification, see \cref{prop:charac_lie_algebra_infitesimal_symmetries} for its characterisation.} The notation emphasises that $\mathfrak{L}_P$ and $\mathfrak{L}_H$ fall into the same class of \cite{bowers2005Classificationthreedimensionalreal} but for different parameters ($2$ and $-1$, respectively), while $\mathfrak{L}_E$ remains apart. When considering complex Lie algebras we would have $\mathfrak{L}_E\cong L\left(3,2,-1\right)$ and thus the three Lie algebras $\mathfrak{L}_Q$ would fall into the same class but in that case we would not be able to distinguish between $\mathfrak{L}_E$ and $\mathfrak{L}_H$ as those complex Lie algebras would be isomorphic. \\
The following proposition gives checkable conditions to identify the Lie algebras $\mathfrak{L}_Q$ among all 3-dimensional Lie algebras. 
\begin{proposition}[Characterisation of $\mathfrak{L}_Q$]
A $3$-dimensional Lie algebra $\mathfrak{L}$ is isomorphic to $\mathfrak{L}_Q$ if and only if $\mathfrak{L}=I\oplus L$, with $I$ an abelian ideal of dimension $2$, and the action $\ad{\ell}{}:I\rightarrow I$ of any $\ell\in L$ on $I$ satisfies 
\begin{enumerate}[label=(\roman*),ref=\textit{(\roman*)}]
    \item either $\ad{\ell}{}$ is diagonalisable with two non-zero real eigenvalues related by $\lambda_1=2\lambda_2$, in which case $\mathfrak{L}\cong\mathfrak{L}_P$.
    \item or $\ad{\ell}{}$ has two purely imaginary eigenvalues, in which case $\mathfrak{L}\cong\mathfrak{L}_E$.
    \item or $\ad{\ell}{}$ is diagonalisable with two non-zero real eigenvalues related by $\lambda_1=-\lambda_2$, in which case $\mathfrak{L}\cong\mathfrak{L}_H$.
\end{enumerate}
\end{proposition}
\begin{remark}\label{rm:lie_ideal}
\tred{Notice that for $\Sigma_E$, resp. $\Sigma_H$, its Lie ideal $\LLLL_0$ is isomorphic to the Lie algebra of infinitesimal symmetries $\LLLL_E$, resp. $\LLLL_H$. This is not the case for the Lie ideal $\LLLL_0$ of $\Sigma_P$ and its symmetry Lie algebra $\LLLL_P$.}
\end{remark}
\begin{proof}
    See \cite{bowers2005Classificationthreedimensionalreal}.
\end{proof}
Before giving our main result, notice that null-forms elliptic and hyperbolic systems $\Sigma_E$ and $\Sigma_H$ do not possess any equilibrium point, i.e. $0\notin\AAA(\xi)$. On the other hand, the drift $f_P=w^2\vec{x}+w\vec{y}$ of the null-form parabolic system $\Sigma_P$ considered locally around $(x_0,y_0,w_0)$ possesses an equilibrium if $w_0=0$ and not if $w_0\neq0$. As a consequence, the system $\Sigma_P$ possesses two non-equivalent local normal forms, which can be represented around $w_0=0$ by 
\begin{multline*}
\Sigma_P^{0}\,:\,\left\{\begin{array}{rl}
    \dot{x} &= w^2  \\
    \dot{y} &= w \\ 
    \dot{w} &= u
\end{array}\right. \\
\textrm{and}\;
\Sigma_P^{1}\,:\,\left\{\begin{array}{rl}
    \dot{x} &=(w+1)^2  \\
    \dot{y} &= w+1 \\ 
    \dot{w} &= u
\end{array}\right.,\;\textrm{respectively}.
\end{multline*}
We now formulate and prove that the Lie algebra $\mathfrak{L}$ of $\Sigma$, being isomorphic to $\mathfrak{L}_Q$, completely characterises those 3-dimensional control-affine systems $\Sigma=(f,g)$, with scalar control, that are locally feedback equivalent to $\Sigma_Q$. Recall that we attach to a control affine system the distributions $\GGG=\distrib{g}$ and $\GGG^1=\distrib{{g},\lb{f}{g}}$.  
\begin{theorem}[$\Sigma_Q$ are characterised by $\mathfrak{L}_Q$] \label{thm:m1_symmetries_characterize_systems}
    Consider the control-affine system $\Sigma:\dot{\xi}=f(\xi)+g(\xi)u$ on a smooth 3-dimensional manifold $\XXX$ with a scalar control, and let $\mathfrak{L}$ be its Lie algebra of infinitesimal symmetries.
    \begin{myenumerate}
        \item[$(i)$]{\label{thm:m1_symmetries_characterize_systems:1} $\Sigma$ is locally feedback equivalent to $\Sigma_E$, around $\xi_0\in\XXX$, if and only if
        \begin{align*}
        \mathfrak{L}\cong\mathfrak{L}_E\quad\textrm{and}\quad I(\xi_0)\oplus\GGG(\xi_0)=T_{\xi_0}\XXX.   
        \end{align*}}
        \item[$(ii)$]{\label{thm:m1_symmetries_characterize_systems:2} $\Sigma$ is locally feedback equivalent to $\Sigma_H$, around $\xi_0\in\XXX$, if and only if
        \begin{align*}
        \mathfrak{L}\cong\mathfrak{L}_H\quad\textrm{and}\quad I(\xi_0)\oplus\GGG(\xi_0)=T_{\xi_0}\XXX.   
        \end{align*}}
        \item[$(iii)$]{\label{thm:m1_symmetries_characterize_systems:3} $\Sigma$ is, locally around $\xi_0\in\XXX$, feedback equivalent to $\Sigma_P^{1}$, around $(x_0,y_0,0)$, if and only if 
        \begin{align*}
            \mathfrak{L}\cong\mathfrak{L}_P,&\quad I(\xi_0)\oplus\GGG(\xi_0)=T_{\xi_0}\XXX,\\
            &\textrm{and}\quad f(\xi_0)\notin\GGG(\xi_0).
        \end{align*}}
        \noindent
        \item[\ref*{thm:m1_symmetries_characterize_systems:3}']{\label{thm:m1_symmetries_characterize_systems:3_p} $\Sigma$ is, locally around $\xi_0\in\XXX$, feedback equivalent to $\Sigma_P^{0}$, around $(x_0,y_0,0)$, if and only if 
        \begin{multline*}
            \mathfrak{L}\cong\mathfrak{L}_P,\quad I(\xi_0)\oplus\GGG(\xi_0)=T_{\xi_0}\XXX,\\
             f(\xi_0)\in\GGG(\xi_0),\quad \textrm{and}\quad \dim\GGG^1(\xi_0)=2.
        \end{multline*}}
    \end{myenumerate}
\end{theorem}
\begin{remark}
     Two systems $\Sigma_E$ (resp. $\Sigma_H$) around any two points $\xi_0$ and $\tilde{\xi}_0$ are locally feedback equivalent to each other. On the other hand, equilibrium points $f(\xi_0)\in\GGG(\xi_0)$ and non-equilibrium points $f(\xi_0)\notin\GGG(\xi_0)$ are distinguished for $\Sigma_P$: the system is feedback equivalent to $\Sigma_P^{0}$ in the former case and to $\Sigma_P^{1}$ in the latter as assert statements \cref{thm:m1_symmetries_characterize_systems:3_p} and \cref{thm:m1_symmetries_characterize_systems:3}, respectively.

    Notice that in statement \cref{thm:m1_symmetries_characterize_systems:3_p} the condition on the pointwise rank of the distribution $\GGG^1$ can be replaced by $(g\wedge \lb{f}{g})(\xi_0) \neq0$, compare this statement with the assumptions of \cref{prop:lp0_equivalent_systems} below.
\end{remark}
\begin{proof}
    We show the sufficiency part of the statements only as the necessity follows immediately from the list of infinitesimal symmetries given by \cref{lem:lie_algebra_symmetries}. For all three cases, the beginning of the proof is the same. 
    
    Consider a control-affine system $\Sigma$, given by vector fields $f$ and $g$, and let three vector fields $v_1$, $v_2$, $v_3$ span the $3$-dimensional Lie algebra $\mathfrak{L}=\vectR{v_1,v_2,v_3}$ of infinitesimal symmetries, which by assumption is isomorphic to $\mathfrak{L}_Q$. We can assume that the abelian ideal of $\mathfrak{L}$ is $I=\vectR{v_1,v_2}$ and that $v_1$,$v_2$,$v_3$ satisfy the commutativity relations of $\mathfrak{L}_Q$ (we express them below for each case separately). Since $I(\xi_0)\oplus \GGG(\xi_0)=T_{\xi_0}\XXX$, it follows that $v_1$, $v_2$, and $g$ are independent, locally around $\xi_0$. We apply a local diffeomorphism $\psi(\xi)=(\tilde{x},\tilde{y},\tilde{w})$, $\psi(\xi_0)=0\in\RR^3$, such that $\tilde{v}_1=\psi_{*}v_1=\vec{\tilde{x}}$, $\tilde{v}_2=\psi_{*}v_2=\vec{\tilde{y}}$, and $\tilde{g}=\psi_*g=g^1\vec{\tilde{x}}+g^2\vec{\tilde{y}}+g^3\vec{\tilde{w}}$, for some smooth functions $g^i$,  satisfying $g^3(0)\neq0$. Replacing $\tilde{g}$ by $\frac{1}{g^3}\tilde{g}$, we may assume that $\tilde{g}=g^1\vec{\tilde{x}}+g^2\vec{\tilde{y}}+\vec{\tilde{w}}$. Then, since $\tilde{v}_1$ and $\tilde{v}_2$ are symmetries of $\GGG= \distrib{g^1\vec{\tilde{x}}+ g^2\vec{\tilde{y}}+ \vec{\tilde{w}}}$, we have $\lb{\tilde{v}_i}{\tilde{g}} \in\GGG$, for $i=1,2$, which implies $g^1=g^1(\tilde{w})$ and $g^2 = g^2(\tilde{w})$. Therefore, we in fact have $\lb{\tilde{v}_1}{\tilde{g}}=\lb{\tilde{v}_2}{\tilde{g}}=0$ and thus there exists a local diffeomorphism $(x,y,w)=\phi(\tilde{x},\tilde{y},\tilde{w})$ such that $\phi_{*}\tilde{v}_1=\vec{x}$, $\phi_{*}\tilde{v}_2=\vec{y}$ and $\phi_{*}\tilde{g}=\vec{w}$. Denote $v_3=v_3^1\vec{x}+v_3^2\vec{y}+v_3^3\vec{w}$, the third infinitesimal symmetry, where $v_3^1=v_3^1(x,y)$ and $v_3^2=v_3^2(x,y)$ since $v_3$ is a symmetry of $\GGG=\distrib{\vec{w}}$. We now separate the case of elliptic, hyperbolic, and parabolic systems.
    \begin{enumerate}[label=\textit{(\roman*)}]
        \item Assume that $\mathfrak{L}\cong\mathfrak{L}_E$ and, since diffeomorphisms do not change the commutation relations, we have $\lb{v_1}{v_3}=-v_2$ and $\lb{v_2}{v_3}=v_1$. Hence, 
        \begin{align*}
            v_3 = (y+c)\vec{x}+(-x+d)\vec{y}+v_3^3(w)\vec{w},
        \end{align*}
        \noindent
        where $c,d$ are real constants. To simplify our computations, we replace $v_3$ by $v_3-c v_1-d v_2\in\mathfrak{L}$, which does not change the commutativity relations, so we can get rid of both constants. Now consider the vector field $f=f^1\vec{x}+f^2\vec{y}+f^3\vec{w}$, then by the fact that $v_1$ and $v_2$ are symmetries of $\Sigma$ we get that $f^1=f^1(w)$ and $f^2=f^2(w)$. Using the symmetry $v_3=y\vec{x}-x\vec{y}+v_3^3(w)\vec{w}$, calculate
        \begin{multline*}
            \lb{v_3}{f}= \left(v_3^3\frac{\partial f^1}{\partial w}-f^2\right)\vec{x}\\
            +\left(v_3^3\frac{\partial f^2}{\partial w}+f^1\right)\vec{y}\mod\GGG,
        \end{multline*}
        \noindent
        which implies the following system of ODEs, where the derivatives are taken with respect to $w$: 
        \begin{align*}
            (Sys_E)\,:\,\left\{\begin{array}{cl}
                v_3^3(f^1)'-f^2 &= 0  \\
                v_3^3(f^2)'+f^1 &= 0  
            \end{array}\right..
        \end{align*}
        \noindent
        By multiplying the first equation of $(Sys_E)$ by $f_1$ and the second by $f_2$ and adding them, we obtain the following relation $v_3^3\left((f^1)^2+(f^2)^2\right)'=0$. From $(Sys_E)$, it is clear that $v_3^3\not\equiv 0$ otherwise we would have $f^1\equiv f^2\equiv 0$ and the Lie algebra of symmetries would be much larger that $\mathfrak{L}$. Observe that the points $w$, where $v_3^3(w)=0$, satisfy also $f^1(w)=f^2(w)=0$. We will show that we always have $\left((f^1)^2+(f^2)^2\right)'=0$. If $v_3^3(0)\neq 0$, then it is clear that the conclusion holds, if $v_3^3(0)=0$ but in an open neighbourhood around $0$ there are no other points where $v_3^3$ vanishes, then it is also clear that the conclusion holds. Finally if $v_3^3(0)=0$ and in any open neighbourhood of zero there are other zeros of $v_3^3$, then consider the open segment $]w_1,w_2[$ between two consecutives zeros. For all $w\in ]w_1,w_2[$ we have $\left((f^1)^2+(f^2)^2\right)'=0$ thus $(f^1)^2+(f^2)^2=R$ ($R\in\RR$) but since for $w=w_1$ and $w=w_2$ we have $f^1(w_i)=f^2(w_i)=0$, for $i=1,2$, therefore, by continuity, $R=0$. But then, we have $\left((f^1)^2+(f^2)^2\right)=0$ implying $f^1=f^2=0$ on $[w_1,w_2]$ and repeating this process for all consecutive roots around $0$ implies that $f^1=f^2=0$ around $0$ but then the Lie algebra $\mathfrak{L}$ of infinitesimal symmetries would be of infinite dimension and this situation is excluded. Finally, we have $\left((f^1)^2+(f^2)^2\right)'=0$ around zero, implying $(f^1)^2+(f^2)^2=R$, $R>0$. Therefore $f^1=\sqrt{R}\cos(w)$ and $f^2=\sqrt{R}\sin(w)$ and, by normalising the coordinates $x$ and $y$, we obtain the form $\Sigma_E$.
        \item The proof is essentially the same so we omit some details. Using the multiplication table $\lb{v_1}{v_3}=v_2$ and $\lb{v_2}{v_3}=v_1$, we obtain 
        \begin{align*}
            v_3=y\vec{x}+x\vec{y}+v_3^3(w)\vec{w},
        \end{align*}
        \noindent
        which, used as a symmetry of $f$, yields the system
        \begin{align*}
            (Sys_H)\,:\,\left\{\begin{array}{cl}
                v_3^3(f^1)'-f^2 &= 0  \\
                v_3^3(f^2)'-f^1 &= 0  
            \end{array}\right..
        \end{align*}
        \noindent
        We derive the relation $v_3^3\left((f^1)^2-(f^2)^2\right)'=0$ and by a similar reasoning we get $\left((f^1)^2-(f^2)^2\right)'=0$ around $0$. Thus, $\left((f^1)^2-(f^2)^2\right)=R$ with $R\in\RR$, but necessarily $R\neq0$ because the Lie algebra of infinitesimal symmetries of $\Sigma$ with $f^1=\pm f^2$ is of infinite dimension. Finally, we have $f^1=\sqrt{|R|}\cosh(w)$ and $f^2=\sqrt{|R|}\sinh(w)$ and using a coordinate normalisation we get $f^1=\cosh(w)$ and $f^2=\sinh(w)$. 
        \item \label{prf:thm:m1_symmetries_characterize_systems:3} Using the multiplication table of $\LLLL_P$ given by $\lb{v_1}{v_3}=2v_1$ and $\lb{v_2}{v_3}=v_2$, which has not been changed by applying diffeomorphisms, we obtain 
        \begin{align*}
            v_3=2x\vec{x}+y\vec{y}+v_3(w)\vec{w},
        \end{align*}
        \noindent
        which is a symmetry of $f$ and thus yields the system 
        \begin{align*}
            (Sys_P)\,:\,\left\{\begin{array}{cl}
                v_3^3(f^1)'-2f^1 &= 0  \\
                v_3^3(f^2)'-f^2 &= 0  
            \end{array}\right.,
        \end{align*}
        \noindent
        recall that $f^1=f^1(w)$ and $f^2=f^2(w)$. We will now distinguish two cases, namely we separate between existence or not of an equilibrium at $w_0=0$. 
        \begin{enumerate}[label=(\alph*),leftmargin=0pt]
            \item Assume that $f(0)\notin\GGG(0)$, that is $(f^1,f^2)(0)\neq(0,0)$ and thus by $(Sys_P)$ we have $((f^1)',(f^2)')(0)\neq (0,0)$. Assume $f^2(0)\neq0$, thus $f^2(w) = c+h(w)$, where $c=f^2(0)$ and $h(0)=0$. Replacing $y$ by $y/c$ we may assume that $f^2(w) =1+h(w)$, where $h'(0)\neq0$, if not the second equation of $({Sys}_P)$ is not satisfied at $0\in\RR$. Set $\hat{w}=h(w)$ and denote the transformed vector fields $\hat{f}$ and $\hat{v}_3$, for which $({Sys}_P)$ implies $\hat{v}_3^3=1+\hat{w}$ and
            \begin{align*}
              \left(\hat{f}^1\right)'(1+\hat{w})&= 2\hat{f}^1.
            \end{align*}
            \noindent
            Solving this equation gives $\hat{f}^1(\hat{w})=c(1+\hat{w})^2$ with $c\in\RR$. But $c$ can not be $0$, otherwise the Lie algebra $\mathfrak{L}$ of infinitesimal symmetries would be of infinite dimension, thus not isomorphic to $\mathfrak{L}_P$. Finally, introducing $\hat{x}=x/c$ we obtain $\Sigma_P^{1}$. If $f^2(0)=0$, then $f^1(0)\neq0$ implying that $(f^1)'(0)\neq0$ and leading to the normalisation $\hat{f}^1(\hat{w})=1+\hat{w}$ giving $\hat{v}_3^3=2(1+\hat{w})$ and $\hat{f}^2=(1+\hat{w})^{1/2}$. This forms is equivalent to $\Sigma_P^{1}$ by the local diffeomorphism $w=(1+\hat{w})^{1/2}-1$, sending $0$ into $0$.
            \item \label{prf:thm:m1_symmetries_characterize_systems:3:b} Assume $f(0)\in\GGG(0)$ and $g\wedge\ad{g}{f}(0)\neq0$. If $(f^2)'(0)\neq0$, take $(\hat{x},\hat{y},\hat{w})=(x,y,f^2(w))$ as a local diffeomorphism around $0\in\RR^3$ that maps $f^1$, $f^2$ and $v_3^3$ into $\hat{f}^1$, $\hat{f}^2$ and $\hat{v}_3^3$, respectively. We have $\hat{f}^2=\hat{w}$, so the system $({Sys}_P)$ implies $\hat{v}_3^3={\hat{w}}$ and 
            \begin{align*}
                \hat{w}\left(\hat{f}^1\right)'&= 2  \hat{f}^1.
            \end{align*}
            \noindent
            Solving this equation gives $\hat{f}^1(\hat{w})=c(\hat{w})^{2}$ with $c\in\RR$. But $c$ can not be $0$, otherwise the Lie algebra $\mathfrak{L}$ of infinitesimal symmetries would be of infinite dimension. Finally, introducing $\hat{x}=x/c$ we obtain $\Sigma_P^{0}$. If $(f^2)'(0)=0$, then $(f^1)'(0)\neq0$ and by applying the local diffeomorphism $(\hat{x},\hat{y},\hat{w})=(x,y,f^1(w))$ we get $\hat{f}^1=\hat{w}$ yielding $\hat{v}_3^3(\hat{w})=2\hat{w}$ and $2\hat{w}(\hat{f}^2)'=\hat{f}^2$. Hence, $|\hat{f}^2|=d|\hat{w}|^{1/2}$ and the only smooth solution, around $\hat{w}=0$, is given by $d=0$ but then the Lie algebra $\mathfrak{L}$ of infinitesimal symmetries would be of infinite dimension contradicting our assumption.
      \end{enumerate}
    \end{enumerate}
\end{proof}
All systems $\Sigma$ that are locally feedback equivalent to either $\Sigma_E$ or $\Sigma_H$ are completely characterised by their symmetry algebras being isomorphic to $\mathfrak{L}_E$ or $\mathfrak{L}_H$, respectively. This is not the case for $\Sigma_P$. Namely, there are systems that have $\mathfrak{L}_P$ as the symmetry algebra although they are not locally feedback equivalent to $\Sigma_P$. 
\begin{proposition}\label{prop:lp0_equivalent_systems}
    Let $\Sigma\,:\,\dot{\xi}=f(\xi)+g(\xi)u$ be a control-affine system on a smooth 3-dimensional manifold $\XXX$ with a scalar control, and let $\mathfrak{L}$ be its Lie algebra of infinitesimal symmetries. Assume $f(\xi_0)\in\GGG(\xi_0)$ and, additionally, that there exists $k\geq1$, the smallest integer such that $g\wedge \adk{g}{k}{f}(\xi_0)\neq0$. Then,
    \begin{align*}
        \mathfrak{L}\cong\mathfrak{L}_P\quad\textrm{and}\quad I(\xi_0)\oplus\GGG(\xi_0)=T_{\xi_0}\XXX
    \end{align*}
    \noindent
    if and only if $\Sigma$ is feedback equivalent to 
    \begin{align*}
        \Sigma_P^{0,k}\,:\,\left\{\begin{array}{cl}
          \dot{z}  &= w^{2k}  \\
          \dot{y}  &= w^k \\ 
          \dot{w}  &= u
        \end{array}\right.
    \end{align*}
    \noindent
    around $(x_0,y_0,0)\in\RR^3$. 
\end{proposition}
Moreover, it is a classical fact that (under the above assumptions) the integer $k$ is an invariant of feedback transformations, hence if $k\neq k'$, ten $\Sigma_P^{0,k}$ and $\Sigma_P^{0,k'}$ are not locally feedback equivalent around $w_0=0$. To be consistent with the notation of the above proposition, the previously considered normal form $\Sigma_P^{0}$ should actually be denoted $\Sigma_P^{0,1}$.
%
\begin{proof}
    We prove the necessity only, as there are no difficulties to show (repeating the arguments of \cref{lem:lie_algebra_symmetries}) that the Lie algebra of infinitesimal symmetries of $\Sigma_P^{0,k}$ is isomorphic to $\mathfrak{L}_P$. We adapt the point \cref{prf:thm:m1_symmetries_characterize_systems:3}-\cref{prf:thm:m1_symmetries_characterize_systems:3:b} of the proof of \cref{thm:m1_symmetries_characterize_systems}.
    Assume that $f(0)\in\GGG(0)$ and that $k\geq1$ is the smallest integer such that $g\wedge\adk{g}{k}{f}\,(0)\neq0$, that is $\frac{\diff^k}{\diff w^k}(f^1,f^2)(0)\neq(0,0)$. If $(f^2)^{(k)}(0)\neq0$, {then by the Taylor expansion we can write $f^2(w)=w^kH(w)$, where $H(0)\neq0$. We can suppose $H(0)>0$, if not, replace $y$ by $-y$, and we apply around $0\in\RR^3$ the local diffeomorphism $(\hat{x},\hat{y},\hat{w})=(x,y,w (H(w))^{1/k})$ that maps $f^1$, $f^2$ and $v_3^3$ into $\hat{f}^1$, $\hat{f}^2$ and $\hat{v}_3^3$, respectively.} 
    We have $\hat{f}^2=\hat{w}^k$, so the system $({Sys}_P)$ implies $\hat{v}_3^3=\hat{w}/k$ and 
    \begin{align*}
        \hat{w}\left(\hat{f}^1\right)'&= 2 k \hat{f}^1.
    \end{align*}
    \noindent
    Solving this singular equation gives $\hat{f}^1(\hat{w})=c\hat{w}^{2k}$ with $c\in\RR$. However, the solution passing through $\hat{w}=0$ is not unique so \`a priori we may have different values of $c$ for $\hat{w}<0$ and $\hat{w}>0$ but the only $C^{\infty}$ solutions are those given by the same value of $c$ (either $c=0$ or $c\neq0$) for any $\hat{w}$. But $c$ can not be $0$, otherwise the Lie algebra $\mathfrak{L}$ of infinitesimal symmetries would be of infinite dimension. Finally, introducing $\hat{x}=x/c$ we obtain the desired form $\Sigma_P^{0,k}$.
    
    If $(f^1)^{(k)}(0)\neq0$, then normalizing $\hat{f}^1=\hat{w}^{k}$ and applying an analogous procedure we deduce that $\hat{f}^2(\hat{w})=c|\hat{w}|^{k/2}$. If $c=0$ then the Lie algebra of infinitesimal symmetries would be of infinite dimension contradicting our assumptions. In all other cases of $k$ and $c$ the solution is not smooth around $\hat{w}=0$ except for $k=2l$ and the same value of $c$ for $\hat{w}<0$ and $\hat{w}>0$. But in the latter case we have $(f^2)^{(l)}(0)\neq0$, with $l<k$, that contradicts the definition of $k$. 
\end{proof}
\begin{remark}
    Statement \cref{thm:m1_symmetries_characterize_systems:3} of \cref{thm:m1_symmetries_characterize_systems} and \cref{prop:lp0_equivalent_systems} describe all systems having $\mathfrak{L}_P$ as the symmetry algebra for which $k$ exists (in particular, all analytic systems). In the $C^{\infty}$ category there are, however, systems for which $k$ does not exist and the symmetry algebra is $\mathfrak{L}_P$. For example, consider 
    \begin{align*}
        \left\{\begin{array}{cl}
          \dot{x}  &= \mathfrak{f}(w)^2   \\
          \dot{y}  &= \mathfrak{f}(w) \\ 
          \dot{w}  &= u
        \end{array}\right.,
    \end{align*}
    \noindent
    with $\mathfrak{f}(w)=\exp\left(-{1}/{w^2}\right)$ and $\mathfrak{f}(0)=0$.
    By a straightforward calculation, its symmetry algebra is, indeed, $\mathfrak{L}_P$ but, obviously, $k$ does not exist at $(x_0,y_0,0)$. 
\end{remark}
{The fact that the case of existence of an equilibrium, i.e. $f(\xi_0)\in\GGG(\xi_0)$, exhibits a richer collection of normal forms can be explained as follows. In that case, the fields of admissible velocities $\AAA=f+\GGG$ becomes at $\xi_0$ a linear (and not an affine) subspace $\AAA(\xi_0)=\GGG(\xi_0)$ of $T_{\xi_0}\XXX$ and the symmetries are less constrained by that structure. Conversely, any given symmetry imposes a bit less rigidity on $f$ if $f(\xi_0)=\GGG(\xi_0)$.}

{In \cite{respondek2002Nonlinearizablesingleinputcontrol} it is proved that any control-affine single-input system, with controllable linear approximation around an equilibrium, has at most two symmetries: the identity  and, in an exceptional "odd" case, a symmetry conjugated to minus identity. All systems studied in the present paper have a nontrivial $3$-dimensional Lie algebra of symmetries but are in a perfect accordance with the result of \cite{respondek2002Nonlinearizablesingleinputcontrol}. Namely, the systems $\Sigma_E$, $\Sigma_H$, and $\Sigma_P^1$ do not have equilibria while $\Sigma_P^0$ and $\Sigma_P^{0,k}$ are considered around an equilibrium but their linear approximations are not controllable. }
\section{Conclusions}

In this paper, we propose a characterisation of null-forms of conic systems $\Sigma_{Q}$ via their Lie algebra of infinitesimal symmetries. We showed that this class of control-affine systems is completely determined by its Lie algebra of symmetries (under some regularity assumptions do the parabolic systems $\Sigma_P$). There are few results of that kind for control-affine systems existing in the literature, with a notable exception of \cite{doubrov2014Geometryrankdistributions}, and more for the control-linear case, see \cite{anderson2011RankdistributionsMonge,doubrov2014modelssubmaximalsymmetric,kruglikov2012LieTheoremRank}. In view of all those results, it would be very interesting to study the following generalisation.
\paragraph{Problem.}
Let $\LLLL$ be a finite-dimensional Lie algebra of vector fields acting transitively on a manifold $\XXX$. Does $\LLLL$ uniquely determine a class of control-affine systems (around a generic point) whose Lie algebra of infinitesimal symmetries is $\LLLL$? \\ 

\noindent
We believe that considering control-affine systems instead of control-linear systems in that problem matters. A reason is that the existence of the drift $f$ gives more rigidity on the symmetries and therefore would constraint the structure of systems possessing those symmetries. 
To start investigating this problem, we could begin with the study of all $2$- and $3$-dimensional Lie algebras (which are well known) and to characterise the control-affine systems that admit those Lie algebras as symmetries. In the thesis of the first author \cite{schmoderer2018Studycontrolsystems}, \cref{thm:m1_symmetries_characterize_systems} is generalised to the class of paraboloid control-affine systems, that is control-affine system (with $m\geq1$ controls and the state space being a $(2m+1)$-dimensional manifold) whose field of admissible velocities is a paraboloid hypersurface of the tangent bundle. 


\printcredits

\bibliographystyle{cas-model2-names}

\bibliography{references/bibliography.bib}

\bio{}
\endbio

\bio{}
\endbio

\end{document}